\documentclass{amsart}[12pt]
\usepackage{amsmath}
\usepackage{amstext}
\usepackage{amssymb}
\usepackage{amsthm}
\usepackage{amscd}

\usepackage{graphicx}

\usepackage{pstcol}
\usepackage{pstricks}
\usepackage{pst-text}
\usepackage{pst-node}
\usepackage{pst-plot}

\theoremstyle{plain}
\newtheorem{thm}{Theorem}[section]
\newtheorem{lem}[thm]{Lemma}
\newtheorem{prop}[thm]{Proposition}
\newtheorem{cor}[thm]{Corollary}

\theoremstyle{definition}
\newtheorem{defi}[thm]{Definition}

\newtheorem{ntn}[thm]{Notation}

\theoremstyle{remark}

 \DeclareMathOperator{\pd}{pd}
\DeclareMathOperator{\Tor}{Tor}
\DeclareMathOperator{\Hred}{\widetilde{H}}
\DeclareMathOperator{\HH}{H} \DeclareMathOperator{\V}{V}
\DeclareMathOperator{\E}{E} \DeclareMathOperator{\link}{link}

\newcommand{\Link}{{\rm{Link}}}
\renewcommand{\b}{\beta}
\newcommand{\e}{\varepsilon}
\newcommand{\AD}{{\Delta^*}}
\newcommand{\linkad}{{\rm{Link}}_{\AD}}
\newcommand{\lc}{ \left\{ }
\newcommand{\rc}{ \right\} }

\begin{document}

\title
[The Betti numbers of forests]
{The Betti numbers of forests}
\author{Sean Jacques}
\address[Jacques]{Department of Pure Mathematics,
University of Sheffield, Hicks Building, Sheffield S3 7RH, United Kingdom\\
{\it Fax number}: 0044-114-222-3769}
\email{PMP00SJ@sheffield.ac.uk}
\author{Mordechai Katzman}
\address[Katzman]{Department of Pure Mathematics,
University of Sheffield, Hicks Building, Sheffield S3 7RH, United Kingdom\\
{\it Fax number}: 0044-114-222-3769}
\email{M.Katzman@sheffield.ac.uk}

\date{\today}

\begin{abstract}
This paper produces a recursive formula of the Betti numbers of certain Stanley-Reisner ideals
(graph ideals associated to forests). This gives a purely combinatorial definition of the
projective dimension of these ideals, which turns out to be a new numerical invariant of forests.
Finally, we propose a possible extension of this invariant to general graphs.
\end{abstract}

\maketitle


\setcounter{section}{-1}
\section{Introduction}

Throughout this paper $K$ will denote a field. For any homogeneous
ideal $I$ of a polynomial ring $R=K[x_1, \dots, x_n]$
there exists a \emph{graded} minimal finite free resolution
\[ 0 \rightarrow \bigoplus_d R(- d)^{\beta_{p d}} \rightarrow \dots
\rightarrow \bigoplus_d R(-d)^{\beta_{1 d}} \rightarrow
R \rightarrow R/I\rightarrow 0\]
of $R/I$, in which $R(-d)$
denotes the graded free module obtained by shifting the degrees of
elements in $R$ by $d$. The numbers $\beta_{i d}$, which we shall
refer to as the $i$th Betti numbers of degree $d$ of $R/I$, are
independent of the choice of graded minimal finite free
resolution. We set $\beta_{0,d}=\delta_{d,0}$ where $\delta$ is
Kronecker's delta, and by convention $\beta_{i,d}=0$ for all $i<0$.
We  also define the total $i$th Betti number of $I$ as
$\beta_i :=\sum \beta_{i d}$.
We refer the reader to chapter 19 of \cite{E}
for an introduction to graded minimal resolutions.


The aim of this paper is to exhibit an interesting combinatorial interpretation of
Betti numbers of graph ideals, which we now define.
Let $G$ be any finite simple graph. We shall always denote the
vertex set of $G$ with $\V(G)$ and its edges with $\E(G)$. Fix an
field $K$ and let $K[\V(G)]$ be the polynomial ring on the vertices of
$G$ with coefficients in $K$. The graph ideal $I(G)$ associated with $G$
is the ideal of $K[\V(G)]$ generated by all degree-2 square-free
monomials $u v$ for which $(u,v)\in \E(G)$. It is not hard to see
that every ideal in a polynomial ring generated by degree-2
square-free monomials is of the form $I(G)$ for some graph $G$.

The quotient $K[\V(G)]/I(G)$ is a always a Stanley-Reisner ring:
define $\Delta(G)$ to be the simplicial complex on the vertices of
$G$ in which a face consists of a set of vertices, no two joined
by an edge. It is easy to see that $K[\V(G)]/I(G)=K[\Delta(G)]$,
the Stanley-Reisner ring associated with $\Delta(G)$.
The simplicial complexes of the form $\Delta(G)$ are characterised by the fact that
their minimal non-faces are one dimensional; these complexes are also known
as \emph{flag complexes.}

Rather than attempting to describe the Betti numbers of graph ideals in terms of the combinatorial
properties of the graph, we shall go the other way around.
The main result in this paper (Theorem \ref{pdTrees}) establishes a new numerical combinatorial
invariant of forests and which is shown to be well defined
by the fact that it coincides with the projective dimension of the ideals associated with forests.
To the best of our knowledge, this is a new invariant of forests.

\section{Hochster's formula}\label{section1}

Recall that for any field $K$ and simplicial complex $\Delta$ the
\emph{Stanley-Reisner ring} $K[\Delta]$ is the quotient of the
polynomial ring in the vertices of $\Delta$ with coefficients in
$K$ by the monomial ideal generated by the product of
vertices not in a face of $\Delta$ (see chapter 5 in \cite{BH} for a good introduction to Stanley-Reisner rings.)

The main tool for investigating Betti numbers of a Stanley-Reisner
ring $K[\Delta]$ is the following theorem giving the Betti numbers as a sum of dimensions
of the reduced homologies of sub-simplicial complexes of $\Delta$.
\begin{thm}[Hochster's Formula (Theorem 5.1 in \cite{H})] \label{thm1.0}
The $i$th Betti number of $K[\Delta]$ of degree $d$ is given by
\[\beta_{i,d} = \sum_{W \subseteq V(\Delta), \#W=d} \dim_K \Hred_{d-i-1}(\Delta_W; K)\]
where $V(\Delta)$ is the set of vertices of $\Delta$ and for any
$W \subseteq V(\Delta)$, $\Delta_W$ denotes the simplicial complex
with vertex set $W$ and whose faces are the faces of $\Delta$
containing only vertices in $W$.
\end{thm}

Notice that when $\Delta=\Delta(G)$ for some graph $G$, we can
rewrite the formula above as
\[
\beta_{i,d} = \sum_{H \subseteq G\mathrm{\ induced\ } \atop{\#\V(H)=d}} \dim_K \Hred_{d-i-1}(\Delta(H); K)
\]

We shall henceforth write $\beta^K_{i,d}(G)$, $\beta^K_{i}(G)$, and $\pd^K(G)$ for the
$i$th Betti number of degree $d$, the $i$th Betti number and the projective dimension of $K[\Delta(G)]$, respectively.
When $K$ is irrelevant or clear from the context, we shall omit the superscript.

As mentioned in the introduction, this paper aims to discover \emph{interesting} combinatorial interpretations of
Betti numbers of graphs.
We would like to mention that these Betti numbers always have \emph{some} combinatorial significance,
as counters (with appropriate weights) of numbers of induced subgraphs in a certain list.
Specifically, if we wanted to interpret $\b^K_{i,d}(G)$ in these terms, we would compile a (finite!)
list $\mathcal{L}_{i,d}$ of all graphs $H$ with $d$ vertices and such that $d_H:=\dim_K \Hred_{d-i-1}(\Delta(H); K)>0$
(this list would depend on $K$; cf.~\cite{K}, for example) and we could write
\[
\b^K_{i,d}(G)=\sum_{H\in \mathcal{L}_{i,d}} n(H) d_H
\]
where $n(H)$ is the number of induced subgraphs  $G$ which are isomorphic to $H$.
For example, for any field $K$
\[
\b^K_{2}(G)=n\left(
\psset{xunit=0.15} \psset{yunit=0.15}
\begin{pspicture}(-0.2,-1)(2,1.2)
\dotnode(0,-1){a1}
\dotnode(2,-1){a3} \dotnode(2,1){a4} \ncline{-}{a1}{a3}
\ncline{-}{a1}{a4}
\end{pspicture}
\right)+
2n\left(
\psset{xunit=0.15} \psset{yunit=0.15}
\begin{pspicture}(-0.2,-1)(2,1.2)
\dotnode(0,-1){a1}
\dotnode(2,-1){a3} \dotnode(2,1){a4} \ncline{-}{a1}{a3}
\ncline{-}{a1}{a4} \ncline{-}{a3}{a4}
\end{pspicture}
\right)+
n\left(
\psset{xunit=0.15} \psset{yunit=0.15}
\begin{pspicture}(-0.2,-1)(2,1.7)
\dotnode(0,-1){a1} \dotnode(0,1){a2} \dotnode(2,-1){a3}
\dotnode(2,1){a4} \ncline{-}{a1}{a2} \ncline{-}{a3}{a4}
\end{pspicture}
\right).
\]

\section{Some elementary properties of Betti numbers}

Recall that the join of two disjoint simplicial complexes $\Delta_1$ and $\Delta_2$, denoted by
$\Delta_1*\Delta_2$, is the simplicial complex with vertices $\V(\Delta_1)\cup \V(\Delta_2)$
and faces
\[
\left\{ F_1 \cup F_2 \,|\, F_1\in \Delta_1, F_2\in \Delta_2 \right\} .
\]
In this section we relate the Betti numbers of the Stanley-Reiner ring
$K[\Delta_1*\Delta_2]\cong K[\Delta_1]\otimes_K K[\Delta_2]$ to those of $K[\Delta_1]$ and $K[\Delta_2]$.

When $G_1$ and $G_2$ are disjoint graphs, $\Delta(G_1\cup G_2)=\Delta(G_1)* \Delta(G_2)$
and the results of this section allow us to deduce that the projective dimension of a graph is additive on its connected
components.

Throughout this section for any simplicial complex $\Delta$,
$\beta_{i,d}^K(\Delta)$  will denote the $i$-th Betti number of degree $d$ of the
Stanley-Reisner ring $K[\Delta]$ as a module over the polynomial ring $K[\V(\Delta)]$ and, similarly,
$\pd^K(\Delta)$ will denote its projective dimension.

\begin{lem}
Let $R=K[x_1, \dots, x_m]$, $S=K[y_1, \dots, y_n]$ be polynomial rings and
let $I\subset R$, $J\subset S$ be homogenous ideals; write $T=R\otimes_K S$.
If $\mathcal{C}$ and $\mathcal{D}$ are minimal graded free resolutions for $R/I$ and $S/J$, respectively, then
$(\mathcal{C} \otimes_R T) \otimes_T  (\mathcal{D}  \otimes_S T) $ is a minimal graded free resolution for $T/(IT+JT)$.
\end{lem}
\begin{proof}
Write $\mathcal{C}' =\mathcal{C} \otimes_R T$ and $\mathcal{D}' =\mathcal{D} \otimes_S T$.
Since $T$ is flat over both $R$ and $S$,
$\mathcal{C}'$ is a $T$-free resolution of $T/IT$
and
$\mathcal{D}'$ is a $T$-free resolution of $T/JT$.

It is not hard to see that
\[
\HH_\bullet\big( \mathcal{C}'  \otimes_T T/JT\big) = \HH_\bullet(\mathcal{C}') \otimes_K T/JT
\]
hence $\Tor_i^T(T/I, T/J)=0$ for all $i>0$.
But we can also compute $\Tor_\bullet^T(T/I, T/J)$ as
$\HH_\bullet\left( \mathcal{C}' \otimes_T  \mathcal{D}' \right)$
(cf.~\cite[\S 6.2]{N})
hence $\mathcal{C}' \otimes_T  \mathcal{D}'$ is a resolution of
$\Tor_0^T(T/I, T/J)=T/(IT+JT)$.

As $\mathcal{C}$ and $\mathcal{D}$ are \emph{minimal} free resolutions, the entries of the maps occurring in
$\mathcal{C}$ and $\mathcal{D}$ (thought of as matrices with entries in $R$ and $S$, respectively) are in the
irrelevant ideals of $R$ and $S$, respectively. This implies that the entries of the maps in
$\mathcal{C}' \otimes_T  \mathcal{D}'$ are in the irrelevant ideal of $T$,  and
so the resolution is minimal.

The standard grading in $R$ and $S$ extends to the standard grading of $T$ and
$\mathcal{C}' \otimes_T \mathcal{D}'$ is easily seen to be graded
with this grading.
\end{proof}

\begin{cor}\label{thm3.0}
Let $\Delta_1$ and $\Delta_2$ be simplicial complexes.
\[
\beta_{i,d}^K(\Delta_1 * \Delta_2) =
\sum_{p+q=i} \sum_{r+s=d} \beta_{p,r}^K (\Delta_1) \beta_{q,s}^K(\Delta_2)
\]
\end{cor}
\begin{proof}
$K[\Delta_1 * \Delta_2]=K[\Delta_1]\otimes_K K[\Delta_2]$.
\end{proof}

\begin{cor}
Let $\Delta_1, \Delta_2$ be disjoint simplicial complexes and let $\Delta=\Delta_1 * \Delta_2 $.
Then $\pd^K(\Delta)=\pd^K(\Delta_1)+\pd^K(\Delta_2)$.
In particular, if $G_1$ and $G_2$ are graphs with disjoint vertices and $G=G_1\cup G_2$, then
$\pd^K(G)=\pd^K(G_1)+\pd^K(G_2)$.
\end{cor}
\begin{proof}
Consider the total Betti numbers
\[
\beta_{i}^K(\Delta_1 * \Delta_2) =
\sum_{p+q=i}  \beta_{p}^K (\Delta_1) \beta_{q}^K(\Delta_2).
\]
We have
\begin{eqnarray*}
\pd^K(\Delta)\geq i & \Leftrightarrow & \beta_{i}^K(\Delta)>0\\
& \Leftrightarrow & \mathrm{there\ exist\ } p,q\geq 0 \mathrm{\ such\ that\ } p+q=i \mathrm{\ and\ } \beta_{p}^K (\Delta_1), \beta_{q}^K(\Delta_2)\neq 0 \\
& \Leftrightarrow & \pd^K(\Delta_1)+\pd^K(\Delta_2) \geq i .\\
\end{eqnarray*}
\end{proof}

\begin{cor}\label{IndependentBetti}
If $\beta_{i,d}^K(\Delta_1)$ and $\beta_{i,d}^K(\Delta_2)$ do not
depend on $K$ for any $i,d$, nor does
$\beta_{i,d}^K(\Delta_1*\Delta_2)$. In particular, if $G_1$ and
$G_2$ are graphs with disjoint vertices and $G=G_1\cup G_2$, then
if $\beta_{i,d}^K(\Delta(G_1))$ and $\beta_{i,d}^K(\Delta(G_2))$
do not depend on $K$ for any $i,d$, nor does
$\beta_{i,d}^K(\Delta(G))$.
\end{cor}

\section{The Eagon-Reiner formula}

The simplicial complexes $\Delta(G)$ do not have an explicit description even for moderately complex graphs $G$ and
as a result, Hochster's formula can not be applied easily to concrete examples.
In \cite{ER} Alexander duality is used to derive a variant of
Hochster's Formula which has the advantage of involving the reduced homologies of
simplicial complexes which are often easier to handle.

Recall that for any simplicial complex
$\Delta$, the \emph{Alexander Dual} of $\Delta$ is the simplicial
complex defined by
\[
\Delta^* := \left\{ F\subseteq \V(\Delta) \,|\, \V(\Delta)\setminus F \notin \Delta \right\} .
\]
The link of a face $F\in \Delta$ is defined as the simplicial
complex
\[
\link_{\Delta} F := \left\{ G\in \Delta \,|\, G\cup F\in \Delta \mathrm{\ and\ } G\cap F=\emptyset \right\} .
\]

\begin{thm}[\cite{ER}]\label{ER}
For all $i\geq 1$, the $\mathbb{N}$-graded Betti numbers of $k[\Delta]$ are given by
\[
\b_{i,d}(K[\Delta])=\sum_{F \in \AD, |F|=n-d} \dim_K \Hred_{i-2}({\Link_{\Delta^*}} F;K) .
\]
\end{thm}

When $\Delta=\Delta(G)$ we write $\Delta^*(G)$ for
$\left(\Delta(G)\right)^*$. Notice that faces of $\Delta^*(G)$ are
the sets of vertices whose complement contains two vertices joined
by an edge in $G$.
For any $F\in \Delta^*(G)$ the simplicial
complex $\link_{\Delta^*} F$ can be described in terms of its maximal faces:
these consist of $\V(G)\setminus (\V(F) \cup \{u,v\})$
for all pair of vertices $u$ and $v$ not in $F$ and
which are connected by an edge in $G$.

\begin{defi}
Let $a_1, \dots , a_s$ be subsets of a finite set $V$. Define
$\e (a_1,\dots,a_s;V)$ to be the simplicial complex which has
vertex set $\bigcup_{i=1}^s (V \setminus a_i)$ and  maximal faces
$V\setminus a_1,\dots, V\setminus a_s$.
(Notice that, with the notation above, the simplicial complex $\e(\{1,2\};\{1,2\})$
has no vertices and is in fact the complex which
has only one face, namely the empty set-- we write this simplicial complex  as $\{ \emptyset\}$ as opposed to $\emptyset$
which has no faces.)
\end{defi}

We can now rephrase the previous remark as follows:
\begin{prop}\label{linkepsilon}
Let $F\in\AD(G)$. Suppose that $e_1,\dots,e_r$ are all the
edges of $G$ which are disjoint from $F$. Then $\linkad F=\e(e_1,\dots,e_r; \V(G) \setminus F)$.
\end{prop}

In the the remainder of this section we establish some homological properties of the simplicial complexes
defined above. The results in this section do not depend on the ground field $K$ and for simplicity we shall write
$\Hred(\bullet)$ for $\Hred(\bullet; K)$.

\begin{lem}\label{lemB}
Let $V$ be a finite set, let
$e_1\dots,e_t,f \in V$ and
write $E_1=\e(e_1\dots,e_t;V)$, $E_2=\e(f;V)$.
We have
$E_1 \cap E_2 =\e(f \cup e_1,\dots,f\cup e_t;V )$.
In particular, if $f \cap \bigcup_{i=1}^t e_i =\emptyset$ then
$E_1 \cap E_2 =\e(e_1,\dots,e_t;V \setminus f )$.
\end{lem}

\begin{proof}
The maximal faces of $E_1 \cap E_2$ are the intersections of the
maximal faces of $E_1$, that is
$V \setminus e_1, \dots V \setminus e_t$, with the maximal face of $E_2$, $V \setminus f$.
These are the sets $V \setminus (e_i \cup f)$ for $i=1,\dots,t$.
Hence we may write
\[
E_1 \cap E_2 =\e(f \cup e_1,\dots,f\cup e_t;V).
\]

Notice that the elements of $f$ are not in any maximal face of $E_1 \cap E_2$ and hence are not vertices of $E_1 \cap E_2$.
If $f \cap \bigcup_{i=1}^t e_i = \emptyset$, we can write
\[
E_1 \cap E_2= \e((e_1 \cup f) \setminus f,\dots,(e_t \cup f) \setminus f;V \setminus f)=\e(e_1,\dots,e_t;V \setminus f ).
\]
\end{proof}

\begin{lem}\label{lemD}
Let $V$ be a finite set, $a\in V$ and $e_1,\dots,e_t\subseteq V\setminus\lc a \rc$.
We have
\[
\Hred_i\big( \e(\{a\},e_1,\dots,e_t;V) \big)=\Hred_{i-1}\big( \e(e_1,\dots,e_t;V \setminus \{a\}) \big)
\]
for all $i$.
\end{lem}

\begin{proof}
Write $E=\e(\{a\},e_1,\dots,e_t;V)$.
Let $\e_1=\e(\{a\};V)$ and let $\e_2=\e(e_1,\dots,e_t;V)$. It is
easily seen that $E=\e_1 \cup \e_2$. The simplicial complex $\e_1$
is in fact just a simplex and so is acyclic.
Also we note that $a \in V \setminus \bigcup_{i=1}^t e_i$ so $a$ is in all maximal faces of $\e_2$ hence
$\e_2$ is also acyclic.
We now make use of the Mayer-Vietoris sequence
\[
\dots \to \Hred_{i}(\e_1) \oplus \Hred_i(\e_2)\to \Hred_i(E)\to \Hred_{i-1}(\e_1 \cap \e_2)\to \dots
\]
which in our case reduces to
\[
\dots \to 0 \to \Hred_i(E) \to \Hred_{i-1}(\e_1 \cap \e_2)\to 0\to \dots
\]
and we obtain $\Hred_i(E) \cong \Hred_{i-1}(\e_1 \cap \e_2)$
for all $i$. By Lemma \ref{lemB}
\[
\e_1 \cap \e_2 =\e(e_1\cup\{a\},\dots,e_t\cup\{a\};V)=\e(e_1,\dots,e_t;V \setminus \{a\})
\]
as $a \notin \bigcup_{i=1}^t e_i$ and we
conclude that
\[
\Hred_i(E) \cong \Hred_{i-1}(\e(e_1,\dots,e_t;V \setminus \{a\} ))
\]
for all $i$.
\end{proof}

\begin{cor}\label{CorE}
Let $V$ be a finite set, let $a_1,\dots,a_s$ be distinct elements of $V$ and let
\[
E=\e(\{a_1\},\dots,\{a_s\},e_1,\dots,e_t;V) .
\]
If
\[
\{a_1,\dots,a_s \} \cap \bigcup_{i=1}^t e_i = \emptyset
\]
then
$\Hred_i(E)=\Hred_{i-s}(\e(e_1,\dots,e_t;V'))$ for all $i$, where $V'=V \setminus \{a_1,\dots,a_s \}$.
\end{cor}

\begin{proof}
Because $a_j \notin (\bigcup_{i=j+1}^s \{a_i\}) \cup
(\bigcup_{i=1}^t e_i)$ for $j=1,\dots, s$ we may repeatedly apply
Lemma (\ref{lemD}) to obtain
\begin{eqnarray*}
\Hred_i(E)&=&\Hred_{i-1}(\e(\{a_2\},\dots,\{a_s\},e_1,\dots,e_t;V
\setminus \{a_1 \}))\\
&=&\Hred_{i-2}(\e(\{a_3\},\dots,\{a_s\},e_1,\dots,e_t;V
\setminus \{a_1,a_2 \}))\\
&=& \dots \\
&=&\Hred_{i-s}(\e(e_1,\dots,e_t;V')).
\end{eqnarray*}
\end{proof}

\section{Betti Numbers of Forests}

Recall that a forest is a graph with no cycles, i.e., a graph whose connected components are trees.
In this section we produce a recursive formula for the Betti numbers of forests in terms of smaller sub-forests.
As a consequence, we obtain an extremely simple recursive formula for the projective dimension of the
Betti numbers of forests.
We thus define a new combinatorial numerical invariant of forests.

We shall refer to the number of neighbours of a vertex $v$ of a graph as the \emph{degree} of $v$.
The crucial property of forests which we will use in this section is the following:
\begin{prop}
Let $T$ be a forest containing a vertex of degree at least two.
There exists a vertex $v$ with neighbours $v_1, \dots, v_n$ where $n\geq 2$ and
$v_1, \dots, v_{n-1}$ and have degree one.
\end{prop}

\begin{proof}
We use induction on the number of vertices of a forest $T$.
Pick a vertex $w$ of degree $1$ and let $T_1=T\setminus \{w\}$.
If $T_1$ has no vertex of degree at least two, take $v$ to be the unique neighbour of $w$ in $T$;
otherwise
the induction hypothesis guarantees the existence of a vertex $u$ in $T_1$ with neighbours $u_1, \dots, u_m$
where $m\geq 2$ and $u_1, \dots, u_{m-1}$ have degree one.
If $(w,u_j)\in \E(T)$ for some $1\leq j\leq m-1$, take $v$ to be $u_j$, otherwise take $v$ to be $u$.
\end{proof}

\begin{ntn}
Henceforth in this section $T$ will denote a forest,
$v$ will be a fixed vertex of $T$ with neighbours $v_1, \dots, v_n$ $(n\geq 2)$ and such that
$v_1, \dots, v_{n-1}$ have degree one.

We will denote with $T'$  the subgraph of $T$ which is obtained by deleting the vertex $v_1$ and
with $T''$  the subgraph of $T$ which is obtained by deleting
the vertices $v,v_1,\dots,v_n$.
Note that $T'$ and $T''$ are both forests.

For a fixed $d>0$ we define the sets
\begin{eqnarray*}
\mathcal{F}_0 & = & \lc F\in \Delta^*(T) \,|\, \V(T'')\subseteq F \rc \\
\mathcal{F}_1 & = & \lc F\in \Delta^*(T) \,|\, v_1 \in F,  |V(T) \setminus F|=d \rc \\
\mathcal{F}_2 & = & \lc F\in \Delta^*(T) \,|\, v_1 \notin F, v\in F,  |V(T) \setminus F|=d \rc \\
\mathcal{F}_3 & = & \lc F\in \Delta^*(T) \,|\, v_1 \notin F, v\notin F,  |V(T) \setminus F|=d \rc \\
\end{eqnarray*}
The sum
\[
 \beta_{i,d}(T)=\sum_{F \in \Delta^*(T): \ |V(T) \setminus F|=d} \dim_K \Hred_{i-2}(\linkad F; K).
\]
giving the graded Betti numbers of $T$ will be split into the sum of
\begin{eqnarray*}
\b_{i,d}^{K}(1) &:=& \sum_{F \in\mathcal{F}_1} \dim_K \Hred_{i-2}(\Link_{\Delta^*(T)} F; K),\\
\b_{i,d}^{K}(2) &:=& \sum_{F \in\mathcal{F}_2} \dim_K \Hred_{i-2}(\Link_{\Delta^*(T)} F; K),\\
\b_{i,d}^{K}(3) &:=& \sum_{F \in\mathcal{F}_3} \dim_K \Hred_{i-2}(\Link_{\Delta^*(T)} F; K).\\
\end{eqnarray*}
\end{ntn}

\begin{lem}\label{lemmaBetta1}
\[
\b_{i,d}^{K}(1)= \beta_{i,d}^K(T').
\]
\end{lem}
\begin{proof}
Consider a face $F$ of $\Delta^* (T)$ with $v_1\in F$.
Then $v_1\notin \Link_{\Delta^* T} F$ and $\Link_{\Delta^* T} F= \Link_{\Delta^* T'} F\setminus \{ v_1 \}$.
Hence
\[
\b_{i,d}^{K}(1)=\sum_{F \in \Delta^*(T'): \ |V(T') \setminus F|=d } \dim_k \Hred_{i-2}({\Link}_{\Delta^*(T')}F;k)=\b_{i,d}(T') .
\]
\end{proof}

\begin{lem}\label{lemmaBetta2}
$\b_{i,d}^{K}(2)=0. $
\end{lem}

\begin{proof}
Using Proposition (\ref{linkepsilon}) we write
$\Link_{\Delta^*(T)} F=\e(e_1,\dots e_r;V)$ for $F \in \Delta^*(T)$, where $e_1,\dots, e_r$
are the edges of $T$ which are disjoint from $F$ and $V= V(T)\setminus F$.

If $F$
includes $v$ but not $v_1$ then $\Link_{\Delta^*(T)} F$ includes $v_1$ but not $v$,
hence the vertex $v$ does not occur in any of the edges $e_1,\dots, e_r$ above and
nor does
$v_1$ belong to these edges since the only edge in $T$
which includes $v_1$ is $\{v,v_1\}$.
Now $v_1\in V$
is in every maximal face of $\Link_{\Delta^*(T)} F$ and so
$\Link_{\Delta^*(T)} F$ is acyclic.
\end{proof}

\begin{lem}\label{lemmaBetta3}
\[
\b_{i,d}^{K}(3) =\sum_{j=0}^{n-1} {{n-1} \choose j} \b_{i-(j+1),d-(j+2)}(T'').
\]
\end{lem}
\begin{proof}
Define a map $\sigma : \Delta^*(T) \rightarrow 2^{\V(T'')}$ by setting
$\sigma(F)=F\setminus \lc v, v_1, \dots, v_n \rc$.
We claim that for any $F\in \Delta^*(T)\setminus \mathcal{F}_0$, $\Link_{\Delta^*(T)} F$ is acyclic or
$\sigma(F)\in \Delta^*(T'')$.
To see this assume that $\sigma(F) \notin \Delta^*(T'')$, i.e., $\V(T)\setminus F$ contains no edge of $T''$.
Since $F\notin  \mathcal{F}_0$, $\Link_{\Delta^*(T)}F$ must contain a vertex of $T''$, and that vertex is in all maximal faces of
$\Link_{\Delta^*(T)} F$ and hence $\Link_{\Delta^*(T)} F$ is acyclic.

Pick any $F\in \mathcal{F}_3$ and write
\[
\Link_{\Delta^*(T)} F= \e( \lc v,v_1 \rc, \lc v,v_{i_1} \rc ,\dots, \lc v,v_{i_j} \rc, e_1, \dots, e_r ;V)
\]
where
$\{v_{i_1},\dots,v_{i_j} \} \subseteq \{v_2,\dots,v_{n} \}$,
$r\geq 0$, $e_1,\dots,e_r$ are  edges which do not feature any of
$v,v_1,\dots,v_n$ and $V=V(T) \setminus F$.
We now show that
\[
\Hred_i\big( \Link_{\Delta^*(T)} F \big) \cong
\Hred_{i-(j+1)}\big( \e(e_1,\dots,e_r; V \setminus \lc v,v_1,v_{i_1},\dots,v_{i_j} \rc )\big)
\]
for all $i$.
Write $E=\Link_{\Delta^*(T)} F$;
we can write $E=\e_1 \cup \e_2$ where
$\e_1=\e(\lc v,v_1 \rc ;V)$ and
\[
\e_2 =\e(\lc v,v_{i_1} \rc ,\dots, \lc v,v_{i_j} \rc, e_1, \dots, e_r ;V).
\]
Now $\e_1$ is a simplex and hence acyclic and
since  $v_1$ is in every maximal face of $\e_2$, $\e_2$ is a cone and, therefore, acyclic.
The corresponding Mayer-Vietoris sequence
\begin{eqnarray*}
\dots \to \Hred_i(\e_1 \cap \e_2)\to \Hred_i(\e_1)\oplus
\Hred_i(\e_2) \to \Hred_i(E) \to \\ \to \Hred_{i-1}(\e_1 \cap
\e_2) \to \Hred_{i-1}(\e_1)  \oplus \Hred_{i-1}(\e_2) \to
\Hred_{i-1}(E) \to \dots
\end{eqnarray*}
reduces to
\begin{eqnarray*}
\dots \to \Hred_i(\e_1 \cap \e_2) \to 0 \to \Hred_i(E) \to
\Hred_{i-1}(\e_1 \cap \e_2) \to 0 \to \Hred_{i-1}(E) \to \dots
\end{eqnarray*}
implying that $\Hred_i(E) \cong \Hred_{i-1}(\e_1 \cap \e_2)$ for all $i$.
The intersection of the simplicial complexes $\e_1$ and $\e_2$ can, by Lemma (\ref{lemB}), be written as
\[
\e_1 \cap \e_2 =\e( \lc v_{i_1} \rc, \dots ,\lc v_{i_j} \rc, e_1,
\dots, e_r; V \setminus \lc v,v_1 \rc).
\]
None of the vertices
$v_{i_1}, \dots, v_{i_j}$  belongs to any of the edges
$e_1,\dots,e_r$ and so by Corollary  (\ref{CorE})
\[
\Hred_i(\e_1 \cap \e_2) \cong \Hred_{i-j}\big(\e(e_1,\dots,e_r;V \setminus \lc v,v_1,v_{i_1},\dots,v_{i_j} \rc) \big)
\]
for all $i$. Putting this together we obtain
\[
\Hred_i(E)=\Hred_{i-1}(\e_1 \cap \e_2)=
\Hred_{i-(j+1)}\big(\e(e_1,\dots,e_r;V \setminus \{v,v_1,v_{i_1},\dots,v_{i_j} \})\big) .
\]

We deduce that, if $F\in \mathcal{F}_3\setminus \mathcal{F}_0$ and $\Link_{\Delta^*(T)} F$ is not acyclic,
we have $\sigma(F)\in \Delta^*(T'')$ and
\[
\Hred_i(\Link_{\Delta^*(T)} F) = \Hred_{i-(j+1)}({\Link}_{\Delta^*(T'')} \sigma(F);  \V(T'')\setminus \sigma(F) ) .
\]
Hence
\[
\sum_{F\in  \mathcal{F}_3\setminus \mathcal{F}_0} \dim_K \Hred_{i-2}\big(\Link_{\Delta^*(T)} F; K\big)=
\sum_{ F\in  \mathcal{F}_3\setminus \mathcal{F}_0\atop{\sigma(F)\in \Delta^*(T^{\prime\prime})}} \Hred_{i-2-(j+1)}\big({\Link}_{\Delta^*(T'')} \sigma(F);\V(T'')\setminus \sigma(F) \big)
\]

Let
$\mathcal{R}=2^{\{ v_2, \dots, v_n \}}$
and for any $\mathbf{\rho} \in \mathcal{R}$
let $U(\mathbf{\rho})= \lc v_2, \dots, v_n\rc \setminus \mathbf{\rho}$.
For any $j\geq 0$ let
\[
\mathcal{L}_j=\big\{ L\in \Delta^*(T'') \,\big|\, \ |\V(T'')\setminus L|=d-(j+2) \big\} .
\]
Notice that for each $L\in \mathcal{L}_j$ and $\mathbf{\rho}\in \mathcal{R}$ with
$| \mathbf{\rho} |=j$ we have
\[
\big| \V(T)\setminus \big(L \cup U(\mathbf{\rho})\big) \big|=
\big| (\V(T'') \setminus L) \cup \{v, v_1\} \cup \mathbf{\rho} \big|=d
\]
and we can now write
\[
\big\{ F\in  \mathcal{F}_3\setminus \mathcal{F}_0 \,\big|\, \ \sigma(F)\in \Delta^*(T^{\prime\prime}) \big\}=
\bigcup_{\mathbf{\rho} \in \mathcal{R}}\
\bigcup_{L\in \mathcal{L}_{|\mathbf{\rho}|}}
L\cup U(\mathbf{\rho})
\]
where all the sets in this union are distinct.

We deduce that
\begin{eqnarray*}
\sum_{F\in  \mathcal{F}_3\setminus \mathcal{F}_0} \dim_K \Hred_{i-2}(\Link_{\Delta^*(T)} F; K)&=&
\sum_{\mathbf{\rho} \in \mathcal{R}}\
\sum_{L\in \mathcal{L}_{|\mathbf{\rho}|}}
\dim_K  \Hred_{i-2}({\Link}_{\Delta^*(T)} (L\cup U(\mathbf{\rho})); K)\\
&=&\sum_{j=0}^{n-1} {n-1 \choose j}
\sum_{L\in \mathcal{L}_j}
\dim_K  \Hred_{i-2-(j+1)}({\Link}_{\Delta^*(T'')} L; K)
\end{eqnarray*}
Notice that for any $L\in \Delta^*(T'')$, $|\V(T'')\setminus L|>0$; using \ref{ER}
we may write the last sum as
\[
\sum_{0\leq j\leq n-1\atop d> j+2} {n-1 \choose j} \beta_{i-(j+1), d-(j+2)} (T'').
\]

Now consider any $F\in \mathcal{F}_3\cap \mathcal{F}_0$;
write
\[
\Link_{\Delta^*(T)} F= \e\big( \lc v,v_1 \rc, \lc v,v_{i_1} \rc ,\dots, \lc v,v_{i_j} \rc ; \V(T)\setminus F\big)
\]
and notice that we must have $d=|\V(T)\setminus F|=j+2$.
It follows from the definition of $\e$ that $\Link_{\Delta^*(T)} F$ is a $(j-1)$-dimensional sphere, so
the homology modules in
\[
\sum_{F\in  \mathcal{F}_3\cap \mathcal{F}_0} \dim_K \Hred_{i-2}(\Link_{\Delta^*(T)} F; K)
\]
vanish unless $i-2=j-1$, i.e., $i=j+1$, in which case they are 1-dimensional.
Thus we can write the sum above as
${n-1 \choose {j}} \delta_{i,j+1}={n-1 \choose {d-2}} \delta_{i,j+1}$ where $\delta_{i,j+1}$ is Kronecker's delta.

We conclude that
\begin{eqnarray*}
\b_{i,d}^{(3)}&=&
\sum_{F\in  \mathcal{F}_3\setminus \mathcal{F}_0} \dim_K \Hred_{i-2}(\Link_{\Delta^*(T)} F; K)+
\sum_{F\in  \mathcal{F}_3\cap \mathcal{F}_0} \dim_K \Hred_{i-2}(\Link_{\Delta^*(T)} F; K)\\
&=&\sum_{0\leq j\leq n-1\atop d>j+2} {n-1 \choose j} \beta_{i-(j+1), d-(j+2)} (T'') + {n-1 \choose {d-2}}\delta_{i,j+1}\\
&=&\sum_{0\leq j\leq n-1} {n-1 \choose j} \beta_{i-(j+1), d-(j+2)} (T'')
\end{eqnarray*}
\end{proof}

\begin{thm}\label{BettiTrees}
\[
\b_{i,d}(T)=\b_{i,d}(T') + \sum_{j=0}^{n-1} {{n-1} \choose {j}} \b_{i-(j+1),d-(j+2)}(T'').
\]
\end{thm}

\begin{proof}
This is immediate from
$\b_{i,d}(T)=\b_{i,d}^{(1)}+\b_{i,d}^{(2)}+\b_{i,d}^{(3)}$
together with Lemmas \ref{lemmaBetta1}, \ref{lemmaBetta2} and \ref{lemmaBetta3}.
\end{proof}

We immediately deduce the following:
\begin{cor}
If $G$ is either
\begin{enumerate}
    \item[(a)] a forest, or
    \item[(b)] a graph whose vertices have degree at most 2,
\end{enumerate}
the Betti numbers of $G$ do not depend on the ground field.
\end{cor}
\begin{proof}
The first case follows directly from Theorem \ref{BettiTrees}.
The connected components of graphs whose vertices have degree at most 2 are paths or cycles.
The Betti numbers of cycles are independent of the ground field (Theorem 7.6.28 in \cite{J}) and the result
follows from Corollary \ref{thm3.0}.
\end{proof}

We now introduce the combinatorial invariant of trees mentioned in the introduction.
Recall that for any graph $G$ and  field $K$ we defined $\pd^K(G)$
as the projective dimension of $K[\Delta(G)]$ as a $K[V(G)]$-module, i.e.,
$\pd^K(G)=\max \lc i\in \mathbb{Z} \,|\, \b^K_{i}(G)>0  \rc$.

\begin{thm}\label{pdTrees}
\[
\pd (T)= \max \lc \pd (T'), \pd (T'')+n \rc.
\]
\end{thm}
\begin{proof}
Theorem \ref{BettiTrees} gives
\begin{eqnarray*}
\b_i(T)&=&\sum_{d \in \mathbb{N}} \b_{i,d}(T)\\
&=& \sum_{d \in \mathbb{N}}\b_{i,d}(T') +
\sum_{d \in \mathbb{N}} \left( \sum_{j=0}^{n-1} {{n-1} \choose j} \b_{i-(j+1),d-(j+2)}(T'')\right).
\end{eqnarray*}
so $\b_i(T)=0$ if and only if $\b_i(T')=0$ and
$\b_{i-1}(T''),\b_{i-2}(T''),\dots,\b_{i-n}(T'')$
are all zero.
But the vanishing of $\b_{i-1}(T''),\b_{i-2}(T''),\dots,\b_{i-n}(T'')$
is equivalent to $\b_{i-n}(T'')=0$.
Hence $\b_i(T)=0$ if and only if $\b_i(T')=0$ and $\b_{i-n}(T'')=0$.
Therefore $i >\pd^K (T)$ if and only if $i>\pd^K (T')$ and
$i-n>\pd^K (T'')$. Hence $\pd^K (T)= \max\{\pd^K (T'),\pd^K (T'')+n\}$.
\end{proof}

It is easy to see that the projective dimension of an isolated vertex is $0$ and that
the projective dimension of an isolated edge is $1$.
Recall also that the invariant $\pd^K(G)$ is additive on connected components (Corollary \ref{IndependentBetti}.)
This together with Theorem \ref{pdTrees} give a complete definition of the projective dimension of all forests $T$.
\emph{The definition of $\pd(T)$ is purely combinatorial.}
If we attempted to define the projective dimension of a forest $T$ using the formula of Theorem \ref{pdTrees}
without being aware of its algebraic significance, it would not be at all clear that the resulting
number is independent of the choice of $v$ which determines the sub-forests $T'$ and $T''$.
We think that the fact that $\pd(T)$ does not depend on these choices is very interesting as there is
no straightforward combinatorial explanation of it.

We now describe another combinatorial invariant of forests which is closely associated with the projective dimension.
\begin{thm}\label{lastBettiTrees}
Let $T$ be a forest. With the notation as in Theorem \ref{pdTrees} we have
\[
\beta_{\pd(T)}(T)=\left\{
\begin{array}{ll}
    \beta_{\pd(T^{\prime})}(T^{\prime}) , & \text{if }  \pd(T^{\prime})> n+ \pd(T^{\prime\prime})\\
    \beta_{\pd(T^{\prime\prime})}(T^{\prime\prime}) , & \text{if }  \pd(T^{\prime})< n+ \pd(T^{\prime\prime})\\
    \beta_{\pd(T^{\prime})}(T^{\prime}) + \beta_{\pd(T^{\prime\prime})}(T^{\prime\prime}) , & \text{if }  \pd(T^{\prime})= n+ \pd(T^{\prime\prime})\\
\end{array}
\right.
\]
\end{thm}
\begin{proof}
This is immediate from  Theorem \ref{BettiTrees} and Theorem \ref{pdTrees}.
\end{proof}
Notice that when $K[\Delta(T)]$ is Cohen-Macaulay $\beta_{\pd(T)}(T)$ is just the
Cohen-Macaulay type of $K[\Delta(T)]$ (cf.~chapter 21 in \cite{E}.)
This invariant is one plus the number of ``ties'' between
$\pd(T^{\prime})$ and $n+ \pd(T^{\prime\prime})$  which occur along any recursive calculation of $\pd(T)$.
We find it striking that this number should be independent of the choices made in this recursive scheme.

\bigskip
At this point it is natural to ask whether Theorem \ref{pdTrees} can be extended to
find a combinatorial interpretation of $\pd^K(G)$ for all graphs $G$.
It is known that $\pd^K(G)$ may depend on the characteristic of the field $K$ (cf.~\cite{K})
and so such an interpretation would have to be ``modular'', in some sense.
But nothing prevents us from defining the following invariant of graphs.
\begin{defi}
Let $G$ be any graph.
Let
\[
p_i(G)=\left\{ T  \,|\, \text{ is a spanning tree in } G \text{ and } \pd(T)=i \right\}
\]
and define the polynomial
$\displaystyle \mathcal{P}_G(x)=\sum_{i\geq 0} |p_i(G)| x^i$.
\end{defi}
The properties of $\mathcal{P}_G(x)$ will be explored in a future paper.

\end{document}